\documentclass[12pt]{amsart}
\usepackage{epsf}
\def\mathbb#1{\mbox{\amss#1}}
\font \amss=msbm10 at 12pt
\usepackage{amsmath}
\usepackage{amssymb}
\usepackage{amsthm}
\newtheorem{teorema}{Theorem}
\newtheorem{predlozhenie}{Proposition}
\newtheorem{lemma}{Lemma}
\newtheorem{defn}{Definition}
\theoremstyle{definition}
\newtheorem{primer}{Example}
\author{E.~Volkov}
\address{Math. Inst. LMU, Theresienstr. 39, 80333 M\"unchen, Germany}
\email{volkov@mathematik.uni-muenchen.de}
\title{Characterization of intrinsically harmonic forms}
\date{June 2007}

\begin{document}

\begin{abstract} Let $M$ be a closed oriented manifold of dimension $n$ 
and $\omega$ a closed $1$-form
on it. We discuss the question whether there exists a Riemannian metric for
which $\omega$ is co-closed. For closed $1$-forms with nondegenerate zeros
the question was answered completely by Calabi in 1969, cf. \cite{Ca}. The 
goal of this paper is to give an answer in the general case, i.e. not making 
any assumptions on the zero set of $\omega$. 
\end{abstract}

\maketitle

\section{Introduction and statement of results}\label{mainres}

In this paper we will be concerned with the characterization
of intrinsically harmonic forms in terms of the 
topological and smooth structure of the underlying manifold. 
A $k$-form on a smooth $n$-manifold is called  
intrinsically harmonic if it is closed and its Hodge-dual with respect to
some Riemannian metric is also closed. So, for a $k$-form to be intrinsically
harmonic it is necessary to be closed and then the desired Riemannian metric 
may or may not exist. The question is: given a closed $k$-form on a smooth
manifold, when is it intrinsically harmonic? Let us give a short historical 
overview. Only the forms of degrees $1$ and $n-1$ have been considered 
seriously. The forms of degrees strictly between $1$ and $n-1$ seem to present 
considerable additional difficulties. 
The following classical theorem of Calabi from 1969
answers the question for forms of degree $1$ with nondegenerate zeros.

\begin{teorema}[Calabi \cite{Ca}]
Let $\omega$ be a closed $1$-form on a 
closed oriented manifold $M$. Assume that all the zeros of $\omega$
are nondegenerate. Then $\omega$ is intrinsically harmonic 
if and only if it is transitive.
\label{Calabi69'}
\end{teorema} 
 
Transitivity will be discussed in detail later in the paper. For now suffice 
it to say that 
a $1$-form is called transitive if there exists a closed transversal
to its kernel foliation through every point which is not a zero of
the form. Note, that one can define the concept of transitivity for
$(n-1)$-forms repeating verbatim the definition above. In 1996 Honda
gave a complete answer to the question for forms with nondegenerate zeros
in degree $n-1$. Having understood the concept of transitivity for both
$1$- and $(n-1)$-forms we can give a unified formulation of the theorems of 
Calabi and Honda. We abbreviate forms with nondegenerate zeros as 
``nondegenerate forms''.

\begin{teorema}[Calabi \cite{Ca}, Honda \cite{Hondath}]
Let $k\in \{1,n-1\}$. For a nondegenerate closed $k$-form $\alpha$ on a 
closed oriented connected $n$-manifold
$M$ to be intrinsically harmonic it is necessary and sufficient that \\
(a) the form $\alpha$ is locally intrinsically harmonic and \\
(b) the form $\alpha$ is transitive.
\label{Calabi69}
\end{teorema} 

We say that a $k$-form is locally intrinsically harmonic, if it becomes
intrinsically harmonic when restricted to a suitable open neighbourhood 
of its zero set. For $1$-forms with nondegenerate zeros transitivity
implies local intrinsic harmonicity and Theorem \ref{Calabi69} simplifies
to Theorem \ref{Calabi69'}. For the discussion of local intrinsic harmonicity
in the case of $(n-1)$-forms see the thesis of Honda \cite{Hondath}. All the 
theorems we have discussed so far assumed that the zeros of the $k$-form
in question are nondegenerate. In 2006 Latschev managed to weaken the 
assumptions on the zero set of the form. 

\begin{teorema}[Latschev \cite{Lat}]
Let $k\in \{1,n-1\}$. Let $\alpha$ be a closed $k$-form on a
closed oriented connected $n$-manifold. Assume that the zero set of $\alpha$
is a Euclidean neighbourhood retract. Then $\alpha$ is intrinsically harmonic
if and only if the following two conditions are satisfied: \\
(a) the form $\alpha$ is locally intrinsically harmonic (if $k=1$ assume in 
addition that the local metric which provides local harmonicity for 
$\alpha$ is real analytic) and \\
(b) the form $\alpha$ is transitive.
\label{Lat}
\end{teorema} 

In this paper we stick to the case of $1$-forms and prove the characterization 
theorem in complete generality, i.e without assuming anything on the local 
structure of the zero set of the $1$-form. We also do not assume any
regularity higher than $C^{\infty}$ for the local metric.

\begin{teorema}
For a closed $1$-form $\omega$ on a closed oriented manifold
$M$ to be intrinsically harmonic it is necessary and sufficient that \\
(a) the form $\omega$ is locally intrinsically harmonic and \\
(b) the form $\omega$ is transitive.
\label{Calabi}
\end{teorema}
For general background we refer to \cite{Fa}.

This work was supported by the DFG grant.

\section{Preliminaries}

We start by discussing the concept of transitivity in a little more detail.
\begin{defn}
A closed $1$-form $\omega$ is called transitive if for any point 
$p\in M\setminus S$ there is a closed strictly $\omega$-positive smooth
(embedded) path $\gamma:S^1\longrightarrow M$ through $p$. 
Here ``strictly $\omega$-positive'' means 
that $\omega(\dot\gamma(t))>0$ for all 
$t\in S^1=\mathbb{R}/\mathbb{Z}$. That is to say that there exists a closed
transversal to the kernel foliation of $\omega$ through every point 
of our manifold which does not lie in the zero set of $\omega$.
\end{defn}

Note that if there exists a smooth, not necessarily embedded, strictly 
$\omega$-positive path through $p$, 
then the path is immersed and if $n=dimM>2$ we can achieve 
embeddedness by a small perturbation. If $n=2$, then we perform the obvious
modifications at double points. In the proof below, however, we get 
embeddedness automatically. To fix conventions from now on 
``$\omega$-positive path'' means ``embedded strictly 
$\omega$-positive path''.   

We recall a classical 
result from dynamical systems --- the Poincar\'e recurrence theorem.
\begin{predlozhenie}
Let $(\Omega, \Sigma, \mu)$ be a probability space. 
Let $\{\phi^t\}_{t \in \mathbb{R}}$ be a measure preserving dynamical system
on it. Assume that $A$ is a $\sigma$-algebra element of positive measure. 
Then for any positive $N$ there exists $n_0$ greater then $N$ such that
$$\mu (A \cap \phi^{n_0}(A))>0.$$
\label{Poincare}
\end{predlozhenie}

To deal with local questions we will need the following
\begin{defn}
For a smooth oriented $n$-dimensional Riemannian manifold $(X,g)$
we define Laplace-Beltrami operator: 
$$\triangle_g:C^{\infty}(X)\longrightarrow \Omega^n(X),$$
which converts a smooth function $f$ into a top degree form 
$$\triangle_gf:=d\star_gdf.$$ 
\label{LB}
\end{defn}

\section{Proof of Theorem \ref{Calabi}}

For necessity assume there exists a 
Riemannian metric $g$ which makes 
$\omega$ harmonic. Condition (a) is obviously 
satisfied. To show Condition (b) we 
%recall a classical 
%result from dynamical systems --- the Poincar\'e-recurrence theorem.
%\begin{predlozhenie}
%Let $(\Omega, \Sigma, \mu)$ be a probability space. 
%Let $\{\phi^t\}_{t \in \mathbb{R}}$ be a measure preserving dynamical system
%on it. Assume that $A$ is a $\sigma$-algebra element of positive measure. 
%Then for any positive $N$ there exists $n_0$ greater then $N$ such that
%$$\mu (A \cap \phi^{n_0}(A))>0.$$
%\label{Poincare}
%\end{predlozhenie}
apply the Poincar\'e recurrence theorem. We set
$\Omega$ to be our manifold 
$M$, the $\sigma$-algebra $\Sigma$ to be the usual borelian $\sigma$-algebra,
and $\mu$ to be the probability measure defined by a distinguished volume 
form $dvol$ on $M$ with total volume equal to one. Furthermore, let the 
vector field $X$ be defined by the following equation: 
$i_Xdvol=\star_{g}\omega$. 
Note that $X$ is transverse to the kernel foliation of $\omega$ outside $S$.
By the 
Cartan formula, we see that $L_Xdvol=0$. Let $\{\phi^t\}_{t \in \mathbb{R}}$
be the flow of $X$ on $M$. In our setting $\{\phi^t\}_{t \in \mathbb{R}}$
becomes a measure preserving dynamical system on $(\Omega, \Sigma, \mu)$.
Let now $p$ be a given point in $M \setminus S$. Let $(\bar {\mathcal{\xi}}, 
\Phi)$ be a bi-foliated closed chart around $p$, 
i.e. $\bar {\mathcal{\xi}}$ is a 
closed subset of $M$, containing an open neighbourhood of $p$ and  
$$\Phi: \bar{\mathcal{\xi}} \longrightarrow B \times I,$$
is a diffeomorphism, 
where $B$ is a closed ball in $\mathbb{R}^{n-1}$ and $I=[0,1]$ 
is a unit time interval. Moreover, under the diffeomorphism 
$\Phi$ flowlines of 
$\{\phi^t\}_{t \in \mathbb{R}}$ correspond to the vertical leaves 
$b \times I$, $b\in B$ and integral submanifolds of the kernel foliation of 
$\omega$ correspond to the horizontal leaves $B \times t$, $t\in I$. In 
further considerations we identify $\bar{\mathcal{\xi}}$ with its image under 
$\Phi$. Since $\bar{\mathcal{\xi}}$ is compact, 
all points of $\bar{\mathcal{\xi}}$ 
will leave it by some time $N$, as we follow the flow 
$\{\phi^t\}_{t \in \mathbb{R}}$. We set $A:=\bar{\mathcal{\xi}}$ and apply 
Proposition \ref{Poincare} with the above choices of 
$\Omega, \Sigma, \mu, A, N$.
This gives us a trajectory of $\{\phi^t\}_{t \in \mathbb{R}}$ which leaves
$\bar{\mathcal{\xi}}$ at some point $(b_1,1)$ and then enters it again for 
the first time at some point $(b_0,0)$. Let us denote the flowline between 
$(b_1,1)$
and $(b_0,0)$ by $\tilde c$. It is clear that except for its end points the 
path $\tilde c$ lies outside $\bar{\mathcal{\xi}}$. Now we close up this 
flowline artificially inside the bifoliated chart $\bar{\mathcal{\xi}}$, 
by connecting $(b_0,0)$ and $(b_1,1)$ with a smooth path $\hat c$ through $p$,
transverse to the horizontal leaves $B \times t$, $t \in I$. Clearly,
this can be done in such a way that the concatenation $c$ of the paths
$\tilde c$ and $\hat c$ is smooth and embedded. So as $c$ is a smooth closed 
$\omega$-positive path and the point $p$ was arbitrary, we have that the 
form $\omega$ is transitive. This is Condition (b).

For sufficiency assume that conditions (a) and (b) hold true. 
Let $U$ be a neighbourhood of $S$ such that $\omega|_{U}$ is co-closed
with respect to some Riemannian metric $g_U$ on $U$. 
It follows from the lemma below that $U$ can be chosen so small that 
the form $\star_{g_U}\omega_{U}$ 
is exact.
\begin{lemma}
Let $(X,g)$ be a smooth oriented $n$-dimensional Riemannian manifold without 
boundary. Let $S$ be a compact zero set of a $1$-form $\gamma$ on $X$ which is
both closed and co-closed. There exists an open neighbourhood $U$ of $S$, 
such that  for any closed $(n-1)$-form $\psi$ on $X$   the restriction 
$\psi|_{U}$ is exact.
\label{algtop}
\end{lemma}
\begin{proof} 
The form $\gamma$ is a solution to a first order linear elliptic equation
\begin{equation}
(d+d^{\star})\gamma=0,
\label{ellipt}
\end{equation}
where $d+d^{\star}=d+\star d\star$ is a Dirac operator on $X$. Locally
\eqref{ellipt} is equivalent to $\triangle_gf=0$, where $f$ is a local 
primitive 
function of $\gamma$.
%and $\triangle_g$ denotes Laplace-Beltrami operator 
%for the metric $g$.
So, we can apply the result by Aronszajn, cf.\cite{Aron}, 
to get that the Dirac operator on $1$-forms possesses 
the strong unique continuation property. Then we apply the 
theorem by C. B\"ar (cf. \cite{Baer}) to find a sequence
$\{L_k\}_{k\in \mathbb{N}}$ of submanifolds of $X$ of codimension at least $2$,
with $S\subset \bigcup_{k\in \mathbb{N}}L_k$. Since every 
submanifold $S_k$ can be countably exhausted by compact ones 
(possibly with boundary), we may without loss of generality assume that each 
$L_k$ is compact, possibly with boundary. Set $Z_k=S\cap L_k$.
Let $dim$ denote the covering dimension of a topological space. Then
$dimZ_k\le n-2$, since $L_k$ is a compact manifold (possibly with boundary)
of dimension at most $n-2$ and $Z_k\subset L_k$. Since 
$S=\bigcup_{k\in \mathbb{N}}Z_k$ and every $Z_k$ is closed in $S$, the 
Countable Sum Theorem  (cf. \cite{Engel} Theorem 7.2.1 on the page 394) 
implies that $dimS\le n-2$. This, in turn, implies that 
$H_{\check{C}ech}^{n-1}(S)=0$. 

Take a sequence 
$\{U_j\}_{j\in \mathbb{N}}$ of open neighbourhoods $U_j$ of $S$ such that
$U_{j+1}\subset U_j$ and $\cap_{j\in \mathbb{N}}U_j=S$ with $U_0=X$. 
The continuity property of \v{C}ech cohomology, cf. \cite{Bred} 
(section 14 ``Continuity'', Theorem 14.4),  
implies that $\varinjlim H_{\check{C}ech}^{n-1}(U_j)=0$, 
but $U_j$ is a manifold, 
hence \v{C}ech cohomology of it is the same as de Rham and finite dimensional. 
So we have that a direct limit of a sequence of 
finite dimensional vector spaces $H_{\check{C}ech}^{n-1}(U_j)$  
vanishes. This implies that for $j$ 
large enough the image of the $0$-th vector space of the sequence 
in the $j$-th one vanishes. In other words if $i:U_j\longrightarrow X$
denotes the obvious inclusion, then $i^{\star}H^{n-1}(X)$ is the trivial 
subspace of $H^{n-1}(U_j)$. Take $U:=U_j$.\end{proof} 
So, we can pick a primitive $(n-2)$-form $\alpha$ on $U$: 
$(d\alpha=\star_{g_U}\omega|_U)$. 
Using transitivity of the form $\omega$, by a 
standard ``thickening of a transversal argument''(see for example \cite{Fa}) 
we obtain that given a point $m \in M \setminus S$, 
there exists an open neighbourhood 
$W_m$ of it, diffeomorphic to 
$S^1\times B$, where $B$ is an open ball in $\mathbb{R}^{n-1}$ centered at the 
origin. Moreover, when restricted to $W_m$, the form $\omega$ is proportional 
to $d\theta$, where $\theta$ denotes the coordinate along the $S^1$ direction.
%More precisely, 
%Let $M$ be a smooth oriented $n$-manifold and $\omega$ --- a $1$-form on it.
%Let $S$ denote the zero set of $\omega$ and 
%let $\gamma:S^1\longrightarrow M$ be a smooth $\omega$-positive path through
%$m \in M \setminus S$, which we have be the transitivity of $\omega$. 
%i.e. 
%$\omega(\dot\gamma(t))>0$ for all $t\in S^1$. 
%Let $W_m\subset M\setminus S$ be a small tubular 
%neighbourhood of $\gamma(S^1)$ in $M$. The neighbourhood $W_m$ 
%is the total space
%of a $D^{n-1}$-bundle $\xi$ over $S^1$, where $D^{n-1}$ is the closed unit 
%disk in $\mathbb{R}^{n-1}$. Every fiber of $\xi$ is a connected component 
%of the intersection of a certain leaf of the kernel foliation of 
%$\omega$ with $W_m$. Since $W_m$ is a total space of a bundle over $S^1$, it 
%can be realized as a mapping torus, i.e. $W_m$ is diffeomorphic to 
%$D^{n-1}\times [0,1]/\sim$, where the equivalence relation is given by  
%$(x,0)\sim (\phi(x),1)$ and $\phi$ is diffeomorphism of $D^{n-1}$. Since
%$W_m$ is orientable, the diffeomorphism $\phi$ is orientation preserving and 
%therefore isotopic to the identity. Therefore, the bundle $\xi$ is trivial,
%i.e. $W_m$ is diffeomorphic to $P=D^{n-1}\times S^1$. For a moment we identify
%$W_m$ with $P$ via this diffeomorphism. Let $x_1,...,x_{n-1}$
%be the coordinates on $D^{n-1}$ and let $\theta$ be the $S^1$ coordinate
%on $P$. In these coordinates the form $\omega|_{W_m}$ writes out as $fd\theta$
%,
%where $f$ is a smooth function on $P$ with $df\wedge d\theta=0$. 
Let 
$\rho:[0,1]\longrightarrow \mathbb{R}$ be a smooth cut-off function:
$\rho|_{[0,1/5]}=1$, $\rho|_{[4/5,1]}=0$. Set 
$\psi_m=\rho(x_1^2+...+x_{n-1}^2)dx_1\wedge...\wedge dx_{n-1}$. Clearly, the 
$(n-1)$-form $\psi_m$ is closed, vanishes 
in a neighbourhood of the boundary of $P$ and the top degree form
$\Theta:=\omega\wedge\psi_m$ satisfies the following properties:
$\Theta$ is nonnegative everywhere and
$\Theta>0$ in some neighbourhood $V_m$ of $\gamma(S^1)$. Vanishing of $\psi_m$ 
near the boundary of $P$ implies that $\psi_m$ vanishes in 
some open neighbourhood $U_m$ of $S$ with $U_m\subset U$. 
This construction almost literally 
follows the one given by Calabi in \cite{Ca}.  

%The 
%co-dimension $1$ foliation on $W_m$ given by the second factor corresponds 
%to the foliation $\mathcal{F}$, a closed transversal $\gamma_j$ through 
%$m$ provided by the transitivity of $\omega$ corresponds to $S^1\times \{0\}$.
%Having this neighbourhood at hands we can construct (see Section \ref{Thom.})
%a closed form $\psi_m$ with the 
%following properties. The form $\psi_m$ vanishes in some neighbourhood 
%$U_m \subset U$ of the set $S$, the top degree form 
%$\omega \wedge \psi_m$ is strictly greater than zero on the closure of 
%some open neighbourhood $V_m\subset W_m$ of $m$ and is nonnegative 
%everywhere. 

Since $M \setminus U$ is compact it can be covered by 
$V_{m_1},...,V_{m_l}$ for some natural number $l$, where 
$m_1,...,m_l \in M \setminus U$. Set $$U_0:=U_{m_1} \cap...\cap U_{m_l},$$
$$V:=V_{m_1} \cup...\cup V_{m_l}$$ and 
$$\psi^{'}:=\Sigma_{i=1}^{l}\psi_{m_i}.$$
Note, that $U_0\subset M\setminus V\subset U$  
and $\psi^{'}|_{U_0}=0$.  

We pause for a moment to summarize what we have. We have
an open neighbourhood $U$ of $S$ with an $(n-2)$-form $\alpha$ 
on $U$ such that $d\alpha=\star_{g_U}\omega$; 
open sets $U_0$ and $V$ with $U_0\subset M\setminus V\subset U$
and an $(n-1)$-form $\psi^{'}$ with $\psi^{'}\wedge \omega$ bounded away form
zero on $V$, nonnegative everywhere and satisfying $$\psi^{'}|_{U_0}=0.$$
This allows us to finish the proof with the standard gluing argument. 
We let $\phi$ be a smooth function with 
$\phi|_{M \setminus V}=1$ and $\phi|_{M\setminus U}=0$.
Such a function $\phi$ exists since both sets $M\setminus V$
and $M\setminus U$ are closed and the first one is contained in the 
complement of the second. Set $\alpha^{''}=\phi \alpha$
and $\psi^{''}=d\alpha^{''}$. Note that $\psi^{''}|_{M\setminus V}=
d\alpha|_{M\setminus V}=\star_{g_{U}}\omega|_{M\setminus V}$.
Consider a closed form $$\psi=K \psi^{'}+\psi^{''}$$
for sufficiently large positive constant $K$. We claim that the form $\psi$ 
has the following properties:
\begin{itemize}
\item [$(i)$]  $(i)\,\, \psi|_{U_0}=\star_{g_U}\omega|_{U_0} $,
\item [$(ii)$] $\omega \wedge \psi>0$ everywhere on $M \setminus S$.
\end{itemize}
Indeed, since $\psi^{'}|_{U_0}=0$, we have that 
$$\psi|_{U_0}=\psi{''}|_{U_0}=\star_{g_{U}}\omega|_{U_0}.$$ 
This shows the first property. 
For the second one consider 
$$\psi_{M\setminus V}=K\psi^{'}|_{M\setminus V}+\psi^{''}|_{M\setminus V}=
K\psi^{'}|_{M\setminus V}+\star_{g_U}|_{M\setminus V},$$
multiplying with $\omega$ gives us 
$$\omega \wedge \psi|_{M\setminus V}=K\omega \wedge \psi^{'}|_{M\setminus V}+ 
\omega \wedge \star_{g_{U}}\omega|_{M\setminus V}.$$ 
The last expression is the sum of 
two nonnegative terms, the second one being strictly positive 
outside $S$. We are left the expression $\omega \wedge \psi$, 
restricted to $V$.
Since $\omega \wedge \psi^{'}|_{V}$ is bounded away from zero, we have that
$$\omega \wedge \psi|_{V}=K \omega \wedge \psi^{'}|_{V}+
\omega \wedge \psi^{''}|_{V}>0$$ for sufficiently large positive constant $K$.

Now, having the form $\psi$ with the properties above we construct the 
desired metric $g$ by gluing. 
Let $\phi_U, \phi_V$ be a partition of unity,
subordinate to the cover $U,V$. Let $g^{''}$ be any metric on $V$,
making $\omega$ and $\psi$ orthogonal to each other. Consider the metric
$\tilde g:=\phi_U g_U+\phi_V g^{''}$ on $M$. 
It makes $\omega$ and $\psi$ orthogonal everywhere on $M$ and 
$\star_{\tilde g}\omega|_{U_0}=\star_{g_U}\omega|_{U_0}=\psi|_{U_0}$. 
Consider the following orthogonal decomposition of the tangent bundle of 
$M \setminus S$:
$$\tilde g=\tilde{g}|_{Ker \omega} \oplus \tilde{g}|_{Ker \psi}.$$
There exists and unique smooth function 
$\tilde{f}: M \setminus S \longrightarrow \mathbb{R}$, such that for the 
metric
$$g:=\tilde{f}\tilde{g}|_{Ker \omega} \oplus \tilde{g}|_{Ker \psi}$$
on $M \setminus S$ we have that 
$\star_g \omega|_{M \setminus S}=\psi|_{M \setminus S}$. Note, that 
$\tilde{f}|_{U_0}=1$, therefore 
$g|_{U_0 \setminus S}=\tilde{g}|_{U_0 \setminus S}$, 
and hence $g$ can be $C^{\infty}$-regularly 
continued across points of $S$ by just setting $g|_S:=\tilde{g}|_S$.
This means that the metric $g$ is well-defined everywhere on $M$. 
The equation $$\star_g \omega=\psi$$ holds on $M \setminus S$, by the 
choice of $\tilde{f}$ and it also holds on $U_0$, by the first property
of the form $\psi$ because $g|_{U_0}=\tilde{g}|_{U_0}=g_U|_{U_0}$. Thus,
since the form $\psi$ is closed we obtain that the form $\omega$
is co-closed with respect to the metric $g$ everywhere on $M$. This completes
the proof of Theorem \ref{Calabi}.

\section{Concluding remarks.}\label{theend}

As we mentioned in the introduction the forms of degrees strictly between
$1$ and $n-1$ present considerable additional difficulties. Indeed,
the simplest case of such a form would be a $2$-form on a $4$-manifold. 
A generic $2$-form on a $4$-manifold does not have any zeros at all. So let 
$\alpha$ be a nowhere zero closed $2$-form on a $4$-manifold. 
To simplify things even further assume that $\alpha$
has constant rank. For dimension reasons we have only two possibilities for 
the rank of $\alpha$ --- $2$ or $4$. In the last case the form $\alpha$ is 
symplectic, and therefore is harmonic for any metric $g$ which is compatible 
with $\alpha$. The question of intrinsic harmonicity is answered trivially
and positively in this case. So the only potentially interesting case is when
$\alpha$ has constant rank $2$. It turns out that this case presents serious
difficulties. As the following example shows,
%The following example was suggested to the author by 
%J.~Latschev. 
%This example shows, that 
transitivity is not sufficient for harmonicity. 

\begin{primer}
Let $M$ be the total space of the nontrivial 
$S^2$-bundle $\xi=(S^2\longrightarrow M\stackrel{\pi}{\longrightarrow} S^2)$ 
over $S^2$. It is easy to see that there exists a section
$s$ of $\xi$ through every point of $M$.
Let $dvol_{S^2}$ be a volume form on the base $S^2$ and set
$\alpha:=\pi^{\star}dvol$. The form $\alpha$ is a closed $2$-form of constant
rank $2$ on
the $4$-dimensional manifold $M$, where the fibers of $\xi$ are the leaves of 
the $2$-dimensional kernel foliation of $\alpha$. Sections of $\xi$ provide
closed $2$-dimensional submanifolds of $M$ to which $\alpha$ restricts as a 
volume form, so $\alpha$ is transitive. But $\alpha$
is not (!) intrinsically harmonic. Assume by contradiction, that there exists
a Riemannian metric $g$ on $M$ such that the form
$\psi:=\star_{g}\alpha$ is closed. The form $\psi$ has constant rank $2$ and
the 
leaves of the kernel foliation of $\psi$ are transverse to those of $\alpha$,
i.e. to the fibers of $\xi$. Take any leaf $\mathcal{L}$ of the kernel 
foliation of $\psi$. The restriction 
$\pi_{\mathcal{L}}:\mathcal{L}\longrightarrow S^2$ is a submersion and
therefore for dimension reasons a covering map. So $\mathcal{L}$ is 
diffeomorphic to $S^2$. So the total space $M$ of $\xi$ admits a foliation
by closed leaves transverse to the fibers with every leaf intersecting 
every fiber exactly once contradicting the nontriviality of $\xi$.
\label{4dim}
\end{primer}
Tautologically one can say that a closed $2$-form of rank $2$ on a $4$-manifold
is intrinsically harmonic if and only if its kernel foliation $Ker\omega$
admits a complementary foliation $\mathcal{F}$ defined by a closed
$2$-form $\psi$. Indeed,
given $g$ which makes $\omega$ harmonic, we can set $\psi:=\star_g\omega$.
Conversely, given a closed $2$-form $\psi$ defining a foliation $\mathcal{F}$
complementary to $Ker\omega$ we can define a Riemannian metric $g$ by requiring
that $Ker\omega$ and $\mathcal{F}$ are orthogonal. Then $\star_g\omega$
is proportional to $\psi$. By rescaling $g$ on $TKer\omega$ we can achieve that
$\star_g\omega$ is equal to $\psi$, which is closed. 
The big problem with this characterization
is that deciding the existence of a complementary foliation defined by a closed
form is just as hard as deciding whether the given form is intrinsically 
harmonic. Finally note that all this is not even beginning to touch the case 
of nonconstant rank.


\begin{thebibliography}{01}

\bibitem{Aron} N.~Aronszajn, A unique continuation theorem for solutions of 
elliptic partial differential equations or inequalities of second order, 
J.~Math. Pures Appl. (9) 36 (1957) pp. 235-249.

\bibitem{Baer} C.~B\"ar, Zero Sets of Solutions to Semilinear Elliptic 
Systems of First Order, Invent. Math., 138 (1999), no. 1, pp. 183-202. 

\bibitem{Bred} G.~Bredon, Sheaf Theory, McGraw-Hill Inc. 1967.
% Section 14 Continuity, Theorem 14.4. 

\bibitem{Ca} E.~Calabi, An intrinsic characterization of harmonic 
one-forms, Global Analysis, Papers in Honour of K.Kodaira, 1969, D.C.~Spencer
and S.~Iyanaga editors, pp. 101-107 Univ.~Tokyo Press, Tokyo.  

\bibitem{Engel} R.~Engelking, General Topology, Helderman Verlag, 1989.
% pp 394 7.2.1. The countable sum theorem

\bibitem{Fa} M.~Farber, Topology of closed One-forms, Mathematical 
Surveys and Monographs vol. 108,  American Mathematical Society 2004.

\bibitem{Hondath} K.~Honda, On Harmonic forms for Generic Metrics, PhD Thesis,
Princeton University, Princeton, 1996.

\bibitem{Lat} J.~Latschev, Closed forms transverse to singular foliations,
Manuscripta Math. 121, 293-315 (2006).
\end{thebibliography}
\end{document}